\documentclass[12pt,reqno]{amsart}

\usepackage{amssymb}

\newcommand{\eps}{\varepsilon}
\renewcommand{\phi}{\varphi}

\newcommand{\Del}{\Delta}

\newcommand{\Z}{{\mathbb Z}}

\newcommand{\cA}{{\mathcal A}}
\newcommand{\cB}{{\mathcal B}}

\newcommand{\oA}{{\bar A}}

\newtheorem{lemma}{Lemma}
\newtheorem{theorem}{Theorem}

\newcommand{\seq}{\subseteq}
\newcommand{\stm}{\setminus}
\newcommand{\est}{\varnothing}

\newcommand{\lpr}{\left(}
\newcommand{\rpr}{\right)}

\newcommand{\longc}{,\dotsc,}
\newcommand{\longp}{+\dotsb+}

\newcommand{\refl}[1]{~\ref{l:#1}}
\newcommand{\reft}[1]{~\ref{t:#1}}
\newcommand{\refs}[1]{~\ref{s:#1}}
\newcommand{\refb}[1]{~\cite{b:#1}}
\newcommand{\refe}[1]{~\eqref{e:#1}}

\author{Vsevolod F. Lev}
\address{Department of Mathematics, The University of Haifa at Oranim,
  Tivon 36006, Israel}
\email{seva@math.haifa.ac.il}

\title[Translation invariance in $\Z/p\Z$]%
  {Translation invariance \\ in groups of prime order}


\makeatletter % to make sure that using the @ symbol in
\newcommand{\subjclassname@NewMSC}%
  {\textup{2010} Mathematics Subject Classification}
\makeatother  % restore the status of @ to the default
\subjclass[NewMSC]{Primary: 11B75; Secondary: 11B25, 11P70}

\begin{document}
\baselineskip 16pt

\begin{abstract}
We prove that there is an absolute constant $c>0$ with the following
property: if $\Z/p\Z$ denotes the group of prime order $p$, and a subset
$A\subset \Z/p\Z$ satisfies $1<|A|<p/2$, then for any positive integer
$m<\min\{c|A|/\ln|A|, \sqrt{p/8}\}$ there are at most $2m$ non-zero elements
$b\in \Z/p\Z$ with $|(A+b)\stm A|\le m$. This (partially) extends onto
prime-order groups the result, established earlier by S.~Konyagin and the
present author for the group of integers.

We notice that if $A\subset \Z/p\Z$ is an arithmetic progression and
$m<|A|<p/2$, then there are exactly $2m$ non-zero elements $b\in \Z/p\Z$ with
$|(A+b)\stm A|\le m$. Furthermore, the bound $c|A|/\ln|A|$ is best possible
up to the value of the constant $c$. On the other hand, it is likely that the
assumption $m<\sqrt{p/8}$ can be dropped or substantially relaxed.
\end{abstract}

\maketitle

\section{Background and motivation}\label{s:bckg}

For a finite subset $A$ and an element $b$ of an additively written
abelian group, let
  $$ \Del_A(b) := |(A+b)\stm A|. $$
If $A$ does not contain cosets of the subgroup, generated by $b$, then the
quantity $\Del_A(b)$ can be interpreted as the smallest number of arithmetic
progressions with difference $b$ into which $A$ can be partitioned. We also
note that $|A|-\Del_A(b)$ is the number of representations of $b$ as a
difference of two elements of $A$; thus, $\Del_A(b)$ measures the
``popularity'' of $b$ as such a difference (with $0$ corresponding to the
largest possible popularity).

The function $\Del_A$ has been considered by a number of authors, the two
earliest appearances in the literature we are aware of being \refb{eh}
and \refb{o}. Evidently, we have $\Del_A(0)=0$; other well-known
properties of this function are as follows:
\begin{itemize}
\item[P1.] $\Del_A(-b)=\Del_A(b)$ for any group element $b$.
\item[P2.] If the underlying group is finite and $\oA$ is the complement
    of $A$, then $\Del_\oA(b)=\Del_A(b)$ for any group element $b$.
\item[P3.] $\Del_A(b_1\longp b_k)\le\Del_A(b_1)\longp\Del_A(b_k)$ for any
    integer $k\ge 1$ and group elements $b_1\longc b_k$.
\item[P4.] Any finite, non-empty subset $B$ of the group contains an
    element $b$ with $\Del_A(b)\ge\lpr1-\frac{|A|}{|B|}\rpr |A|$.
\end{itemize}
The interested reader can find the proofs in \cite{b:eh,b:o,b:hls} or work
them out as an easy exercise. We confine ourselves to the remark that the
last property follows by averaging over all elements of $B$.

The basic problem arising in connection with the function $\Del_A$ is to
show that it does not attain ``too many'' small values; that is, every
set $B$ contains an element $b$ with $\Del_A(b)$ large, with the precise
meaning of ``large'' determined by the size of $B$. Accordingly, we let
  $$ \mu_A(B) := \max_{b\in B} \Del_A(b). $$
Property P4 readily yields the simple lower-bound estimate
\begin{equation}\label{e:mu-lower-simple}
  \mu_A(B) \ge \lpr 1-\frac{|A|}{|B|} \rpr |A|;
\end{equation}
however, this estimate is far from sharp, and insufficient for most
applications.

Notice, that if $d$ is a group element of sufficiently large order, $A$ is an
arithmetic progression with difference $d$, and $B=\{d,2d\longc md\}$ with
$m=|B|\le|A|$, then $\mu_A(B)=|B|$. Thus,
\begin{equation}\label{e:mu-ge-|B|}
  \mu_A(B) \ge |B|
\end{equation}
is the best lower-bound estimate one can hope to prove under the assumption
$B\cap(-B)=\est$ (cf.~Property~P1). In view of the trivial inequality
$\mu_A(B)\le |A|$, a necessary condition for \refe{mu-ge-|B|} to hold is
$|B|\le|A|$, but this may not be enough to require: say, an example presented
in \refb{kl} shows that \refe{mu-ge-|B|} fails in general for the group of
integers, unless $|B|<c|A|/\ln|A|$ with a sufficiently small absolute
constant $c$. As shown in \refb{kl}, this last assumption already suffices.

\begin{theorem}[{\cite[Theorem~1]{b:kl}}]\label{t:kl}
There is an absolute constant $c>0$ such that if $A$ is a finite set of
integers with $|A|>1$, and $B$ is a finite set of positive integers
satisfying $|B|<c|A|/\ln|A|$, then $\mu_A(B)\ge|B|$.
\end{theorem}

\section{The main result}

It is natural to expect that an analogue of Theorem~\reft{kl} remains valid
for groups of prime order, particularly since the arithmetic progression case
is ``worst in average'' for these groups: namely, it is easy to derive from
\cite[Theorem~1]{b:l2} that for all sets $A$ and $B$ of given fixed size in
such a group, satisfying $B\cap(-B)=\est$, the sum $\sum_{b\in B}\Del_A(b)$
is minimized when $A$ is an arithmetic progression, and $B=\{d,2d\longc
md\}$, where $m$ is a positive integer and $d$ is the difference of the
progression. The goal of this note is to establish the corresponding
supremum-norm result.

Throughout, we denote by $\Z$ the group of integers, and by $\Z/p\Z$ with $p$
prime the group of order $p$.

\begin{theorem}\label{t:main}
There exists an absolute constant $c>0$ with the following property: if $p$
is a prime and the sets $A,B\subset \Z/p\Z$ satisfy
$1<|A|<p/2,\,B\cap(-B)=\est$, and $|B|<\min\{c|A|/\ln|A|,\sqrt{p/8}\}$, then
$\mu_A(B)\ge|B|$.
\end{theorem}

As Property P2 shows, the assumption $|A|<p/2$ of Theorem~\reft{main} does
not restrict its generality. In contrast, the assumption $|B|<\sqrt{p/8}$
seems to be an artifact of the method and it is quite possible that the
assertion of Theorem~\reft{main} remains valid if this assumption is
substantially relaxed or dropped altogether.

We notice that Theorem~\reft{main} is formally stronger than
Theorem~\reft{kl}. However, the proof of the former theorem (presented in
Section~\refs{proof}) relies on the latter one, used ``as a black box''. The
proof also employs a rectification result of Freiman, and elements of the
argument used in \refb{kl} to prove Theorem~\reft{kl}, in a somewhat modified
form.

The rest of this paper is divided into three parts: having prepared the
ground in the next section, we prove Theorem~\reft{main} in
Section~\refs{proof}, and present an application to the problem of estimating
the size of a restricted sumset in the last section.

%. In the next section we collect some auxiliary results, needed in the course
%of the proof of Theorem \reft{main}; the theorem In Section \refs{rsum} we
%apply the estimates of the function $\mu_A$ to derive two results on the
%smallest possible size of a restricted sumset.

%The proof of Theorem~\reft{main}, presented in Section \refs{proof}, combines
%a rectification result of Freiman, Theorem~\reft{kl} used ``as a black box'',
%and elements of the argument employed in \refb{kl} to prove the latter
%theorem (in a somewhat modified form). Auxiliary results needed in the course
%of the proof are collected in the next section.

\section{The toolbox}

In this section we collect some auxiliary results, needed in the course
of the proof of Theorem \reft{main}.

Given a subset $B$ of an abelian group and an integer $h\ge 1$, by $hB$ we
denote the \emph{$h$-fold sumset of $B$}:
  $$ hB := \{b_1\longp b_h \colon b_1\longc b_h\in B \}. $$

Our first lemma is an immediate consequence of Property P3.
\begin{lemma}\label{l:subadditivity}
For any integer $h\ge 1$ and finite subsets $A$ and $B$ of an abelian group
we have
  $$ \mu_A(hB) \le h\mu_A(B). $$
\end{lemma}

The following lemma of Hamidoune, Llad\'o, and Serra gives an estimate which,
looking deceptively similar to \refe{mu-lower-simple}, for $B$ small is
actually rather sharp. We quote below a slightly simplified version, which is
marginally weaker than the original result.

\begin{lemma}[{\cite[Lemma~3.1]{b:hls}}]\label{l:hls}
Suppose that $A$ and $B$ are non-empty subsets of a finite cyclic group such
that $B\cap(-B)=\est$ and the size of $A$ is at most half the size of the
group. If every element of $B$ generates the group, then
  $$ \mu_A(B) > \lpr 1-\frac{|B|}{|A|}\rpr |B|. $$
\end{lemma}

Yet another ingredient of our argument is a rectification theorem due to
Freiman.
\begin{theorem}[{\cite[Theorem~2.11]{b:n}}]\label{t:f}
Let $p$ be a prime and suppose that $B\subset\Z/p\Z$ is a subset with
$|B|<p/35$. If $|2B|\le 2.4|B|-3$, then $B$ is contained in an arithmetic
progression with at most $|2B|-|B|+1$ terms.
\end{theorem}

Finally, we need a lemma showing that if $B$ is a dense set of integers, then
the difference set
  $$ B-B := \{ b'-b'' \colon b',b''\in B \} $$
contains a long block of consecutive integers.
\begin{lemma}[\protect{\cite[Lemma~3]{b:l}}]\label{l:difset}
Let $B$ be a finite, non-empty set of integers. If
 $\max B-\min B<\frac{2k-1}k\,|B|-1$ with an integer $k\ge 2$, then $B-B$
contains all integers from the interval $(-|B|/(k-1),|B|/(k-1))$.
\end{lemma}

\section{Proof of Theorem~\reft{main}}\label{s:proof}

For real $u<v$ and prime $p$, by $\phi_p$ we denote the canonical
homomorphism from $\Z$ onto $\Z/p\Z$, and by $[u,v]_p$ the image of the set
$[u,v]\cap\Z$ under $\phi_p$. In a similar way we define $[u,v)_p$ and
$(u,v)_p$.

We begin with the important particular case where $B$ is a block of
consecutive group elements, starting from $1$. Thus, we assume that $p$ is a
prime, $A\subset \Z/p\Z$ satisfies $1<|A|<p/2$, and
$m<\min\{c|A|/\ln|A|,\sqrt{p/8}\}$ is a positive integer (where $c$ is the
constant of Theorem~\reft{kl}), and show that, letting then $B:=[1,m]_p$, we
have $\mu_A(B)\ge m$.

Suppose, for a contradiction, that $\mu_A(B)<m$. Since $A$ is a union of
$\Del_A(1)$ blocks of consecutive elements of $\Z/p\Z$, so is its complement
$\oA:=(\Z/p\Z)\stm A$, and we choose integers $u<v$ such that
$[u,v)_p\seq\oA$ and
\begin{equation}\label{e:big-gap}
  v-u \ge \frac{|\oA|}{\Del_A(1)} > \frac{p}{2m} > m.
\end{equation}
Rectifying the circle, we identify $A$ with a set of integers
$\cA\seq[v,u+p)$, and $B$ with the set $\cB:=[1,m]\cap\Z$. Inequality
\refe{big-gap} shows that an arithmetic progression in $\Z/p\Z$ with
difference $d\in[1,m]_p$ cannot ``jump over'' the block $[u,v)_p$; hence,
$\mu_A(B)=\mu_\cA(\cB)$. On the other hand, we have $\mu_\cA(\cB)\ge|\cB|=m$
by Theorem~\reft{kl}. It follows that $\mu_A(B)\ge m$, the contradicting
sought.

We notice that so far instead of $m<\sqrt{p/8}$ we have only used the weaker
inequality
\begin{equation}\label{e:sqrt}
  m<\sqrt{p/2};
\end{equation}
this observation is used below in the proof.

Having finished with the case where $B$ consists of consecutive elements of
$\Z/p\Z$, we now address the general situation. Suppose, therefore, that
$A,B\seq \Z/p\Z$ satisfy the assumptions of the theorem and, again, assume
that $\mu_A(B)<|B|$.

For a subset $S$ of an abelian group we write $S^\pm:=S\cup\{0\}\cup(-S)$;
thus, by Property~P1, we have $\mu_A(S^\pm)=\mu_A(S)$ for any finite subset
$A$ of the group, and if $S\cap(-S)=\est$, then $|S^\pm|=2|S|+1$.

If $|2B^\pm|\ge \frac13\,|A|+1$, then by the well-known Cauchy-Davenport
inequality (see, for instance, \cite[Theorem~2.2]{b:n}), we have
$|12B^\pm|>2|A|$. Thus, using Lemma~\refl{subadditivity} and estimate
\refe{mu-lower-simple}, and assuming that $c$ is sufficiently small, we
conclude that
  $$ \mu_A(B) = \mu_A(B^\pm) \ge \frac1{12}\,\mu_A(12B^\pm)
                                               > \frac1{24}\,|A| \ge |B|, $$
a contradiction; accordingly, we assume
  $$ |2B^\pm|<\frac13\,|A|+1. $$
Let $C:=(2B^\pm)\cap[1,p/2)_p$. Observing that
$|C|=(|2B^\pm|-1)/2<\frac16\,|A|$, by Lemmas~\refl{hls}
and~\refl{subadditivity} and the assumption $\mu_A(B)<|B|$ we get
  $$ \frac56\,|C| < \mu_A(C) = \mu_A(2B^\pm) \le 2\mu_A(B^\pm)
                        = 2\mu_A(B) \le 2(|B|-1) = |B^\pm|-3; $$
hence,
\begin{equation}\label{e:3Bpmdense}
  |2B^\pm| = 2|C| + 1 < \frac{12}5\,|B^\pm| - \frac{31}{5} < 2.4|B^\pm|-3.
\end{equation}
We now apply Theorem~\reft{f} to derive that the set $B^\pm$ is contained in
an arithmetic progression with at most
$|2B^\pm|-|B^\pm|+1<\frac13\,|A|<p/2+1$ terms. Taking into account that $0\in
B^\pm$ and dilating $A$ and $B$ suitably, we assume without loss of
generality that $B^\pm\seq(-p/4,p/4)_p$ and $B^\pm$ is actually contained in
a block of at most $|2B^\pm|-|B^\pm|+1$ consecutive elements of $\Z/p\Z$.

Let $\cB\seq[1,p/4)$ be the set of integers such that
$B^\pm=\phi_p(\cB^\pm)$, and write $l:=\max(\cB^\pm)-\min(\cB^\pm)$. From
\refe{3Bpmdense} we conclude that
  $$ l \le |2B^\pm|-|B^\pm| < \frac32\, |B^\pm| - 1
                                                 = \frac32\, |\cB^\pm| - 1. $$
Therefore, by Lemma~\refl{difset} (applied with $k=2$) we have
  $$ [1,|\cB^\pm|-1] \seq \cB^\pm-\cB^\pm = 2\cB^\pm, $$
whence
  $$ [1,|B^\pm|-1]_p \seq 2B^\pm. $$

Recalling that the result is already established for the consecutive residues
case, and observing that $|B^\pm|-1=2|B|<\sqrt{p/2}$ (to be compared with
\refe{sqrt}), we obtain
  $$ \mu_A(2B^\pm) \ge \mu_A([1,|B^\pm|-1]_p) \ge |B^\pm|-1 = 2|B|. $$
Using now Lemma~\refl{subadditivity} we get
  $$ 2\mu_A(B) = 2\mu_A(B^\pm) \ge \mu_A(2B^\pm) \ge 2|B|, $$
a contradiction completing the proof of Theorem~\reft{main}.

\section{An application: restricted sumsets in abelian groups}\label{s:rsum}

\newcommand{\tauplus}{\stackrel{\tau}{+}}

Given two subsets $A$ and $B$ of an abelian group and a mapping
 $\tau\colon B\to A$, let
  $$ A\tauplus B := \{a+b\colon a\in A,\,b\in B,\,a\neq\tau(b) \}. $$
Restricted sumsets of this form, generalizing in a natural way the
``classical'' restricted sumset $\{a+b\colon a\in A,\,b\in B,\,a\neq b\}$,
were studied, for instance, in \refb{l3}. Since
  $$ |(A+b_1)\cup(A+b_2)| = |A|+|(A+b_1-b_2)\stm A| $$
for any $b_1,b_2\in B$, we have
\begin{align*}
  |A+B| &\ge |A| + \mu_A(B-B)
\intertext{and, furthermore,}
  |A\tauplus B| &\ge |A| + \mu_A(B-B) - 2;
\end{align*}
hence, lower-bound estimates for $\mu_A(B-B)$ translate immediately into
estimates for the cardinalities of the sumset $A+B$ and the restricted sumset
$A\tauplus B$. Here we confine ourselves to stating three corollaries of
estimate \refe{mu-lower-simple}, Lemma~\refl{hls}, and Theorem~\reft{main},
respectively.

\begin{theorem}
Suppose that $A$ and $B$ are finite subsets of an abelian group. If for some
real $\eps>0$ we have $|B|\le(1-\eps)|A|$ and $|B-B|\ge\eps^{-1}|A|$, then
\begin{align*}
  |A+B| &\ge |A|+|B|
\intertext{and}
  |A\tauplus B| &\ge |A|+|B|-2
\end{align*}
for any mapping $\tau\colon B\to A$.
\end{theorem}

\begin{theorem}
Suppose that $p$ is a prime and $A,B\seq\Z/p\Z$ are non-empty. If $|A|<p/2$
and $|B|<\sqrt{|A|}+1$, then for any mapping
 $\tau\colon B\to A$ we have
  $$ |A\tauplus B| \ge |A|+|B|-3. $$
\end{theorem}

For the proof just notice that if $2\le|B|\le(p+1)/2$, then by the
Cauchy-Davenport inequality there exists a subset $C\seq B-B$ with
$C\cap(-C)=\est$ and $|C|=|B|-1$, whence, in view of Lemma \refl{hls},
  $$ \mu_A(B-B) \ge \mu_A(C) \ge \left(1-\frac{|C|}{|A|}\right) |C|
       = |B| - 1 - \frac{(|B|-1)^2}{|A|} > |B| - 2. $$

\begin{theorem}
Suppose that $p$ is a prime and $A,B\seq\Z/p\Z$. If $1<|A|<p/2$ and
$0<|B|<\min\{\sqrt{p/8},\,c|A|/\ln|A|\}$, where $c$ is a positive absolute
constant, then for any mapping $\tau\colon B\to A$ we have
  $$ |A\tauplus B| \ge |A|+|B|-3. $$
\end{theorem}

In connection with the last two theorems we notice that a construction
presented in \refb{l3} shows that for (non-empty) subsets $A,B\seq\Z/p\Z$ and
a mapping $\tau\colon B\to A$, the estimate $|A\tauplus B|\ge|A|+|B|-3$ may
fail in general, even if the right-hand side is substantially smaller than
$p$. A question raised in \refb{l3} and ramaining open till now is whether
this estimate holds true under the additional assumption that $\tau$ is
injective and $|A|+|B|\le p$.

\bigskip
\end{document}